\documentclass[11pt, a4paper]{article}  

\author{Frank Schuhmacher}

\title{Convergent noncommutative power series algebras and their 
finite dimensional representations}

\usepackage{mathenv}

\usepackage[intlimits]{amsmath}

\usepackage[T1]{fontenc}

\usepackage{amssymb}

\usepackage{amsthm}

\usepackage{amscd}

\usepackage[all]{xypic}

\usepackage{makeidx}

\newcounter{punkt}

\theoremstyle{definition}

\newtheorem{defi}{Definition}[section]

\newtheorem{beisp}[defi]{Example}

\theoremstyle{theorem}

\newtheorem{satz}[defi]{Theorem}

\newtheorem{lemma}[defi]{Lemma}

\newtheorem{kor}[defi]{Corollary}

\newtheorem{prop}[defi]{Proposition}

\newtheorem{observe}[defi]{Observation}

\newcommand{\nach}{\longrightarrow}

\newcommand{\sub}{\subseteq}

\newcommand{\eins}{\text{\textbf{1}}}
\newcommand{\sgn}{\operatorname{sgn}}

\newcommand{\MM}{ \mathbb{M} }
\newcommand{\RR}{ \mathbb{R} }

\newcommand{\ZZ}{\mathbb{Z}}

\newcommand{\CC}{\mathbb{C}}

\newcommand{\C}{\mathcal{C}}

\renewcommand{\O}{\mathcal{O}}

\newcommand{\M}{\mathcal{M}}

\newcommand{\ppp}{\cdot\ldots\cdot}

\newcommand{\sC}{\mathcal{C}}
\newcommand{\K}{\mathcal{K}}

\newcommand{\m}{\mathfrak{m}}

\newcommand{\spec}{\operatorname{Spec}}

\newcommand{\Hom}{\operatorname{Hom}}

\newcommand{\Kern}{\operatorname{Kern} }

\newcommand{\ot}{\otimes}

\newcommand{\sig}{\sigma}

\newcommand{\ddd}{,\ldots,}

\newcommand{\kkk}{+\ldots+}

\newcommand{\bbb}{|\ldots|}
\newcommand{\ab}{\operatorname{ab}}
\newcommand{\abs}{\operatorname{abs}}
\newcommand{\unab}{\operatorname{unab}}
\newcommand{\op}{\operatorname{op}}
\newcommand{\ec}{\hat{e}}

\newcommand{\alphac}{\hat{\alpha}}
\newcommand{\otc}{\hat{\otimes}}

\begin{document}
\maketitle
\begin{abstract}
We define ``convergence'' for 
noncommutative power series and construct 
two topologies on the algebra of power series,
convergent with respect to a positive radius. We 
indicate all finite dimensional continuous
representations of this algebra and
prove completeness for both topologies. 
\end{abstract}
\tableofcontents

\section*{Introduction}

In \cite{ncaspaces}, the author constructed
a sheaf of associative $\CC$-algebras on 
the space $\CC^n$ with the sheaf of holomorphic functions
as abelianization and ``convergent'' noncommutative
power series algebras as local rings.
This construction led naturally to a definition of
``noncommutative analytic spaces'' as a generalization of
(real or complex) analytic spaces. Examples for
noncommutative analytic spaces
are analytic super-manifolds
and noncommutative links.\\

This paper starts with a more refined definition of
``convergence'' for noncommutative power series. 
A series, convergent in the sense
of \cite{ncaspaces}, is called ``absolutely convergent'' here
and in the future.
To explain the refined definition,
consider the space $\MM$ of complex $k\times k$-matrices.
The star shaped, open subsets 
$$P^{\MM}(r)$$
of $\MM^n$ of all $n$-tuples of matrices with
``joint spectral radius'' smaller than a positive $r$,
play a crucial role. Roughly speaking,
a noncommutative power series $f$ is convergent with
respect to $r$, if all the series 
$f(M)$ obtained by plugging in elements $M$ of $P^{\MM}(r)$ in $f$,
converge, for arbitrary matrix length $k$.
(In the subsequent paper \cite{spaces},
we show that this notion of convergence can be used
to define a structure sheaf $\O$ on $\CC^n$ which is slightly richer
than the sheaf defined in \cite{ncaspaces}.)
We prove (Theorem~\ref{frechet})
that the algebra
$$\CC\{x\}_r=\CC\{x_1|...|x_n\}_r$$
of powers series, convergent with respect to
a given $r>0$, is Fr\^echet
in the ``topology of compact convergence'', where
``compact'' refers to subsets of $P^{\MM}(r)$.
We show (Theorem~\ref{volls}) that
$\CC\{x\}_r$
is complete with respect to the filtration
by its ``ideals of matrix identities''.
An element $f\in\CC\{x\}_r$ is a ``matrix identity'' for
the matrix algebra $\MM$, if
$f(M)=0$, for each $M\in P^{\MM}(r)$.
On the geometric side, 
we prove (Theorem~\ref{homeom}) that the continuous $k$-dimensional
representation of $\CC\{x\}_r$ are just those of 
the form $f\mapsto f(M)$, where $M$ is
in $P^{\MM}(r)$. \\

\textbf{Conventions:} The letter $\CC$ denotes weather the fiel of real
or the field of complex numbers. We will indicate explicitly if a statement 
only holds for the complex numbers.
An \textit{algebra} is  an associative, unitary $\CC$-algebra.
\textit{Modules} and \textit{ideals} are two-sided. 
Left-modules only appear under the name ``representations''.
An element $a$
of an algebra $A$ is a \textit{unit}, if there exists an element 
$b\in A$ such that $ab=1=ba$. 
A \textit{maximal ideal} is maximal in the set of (two-sided) ideals.
An algebra $A$ with a single maximal ideal $\m$ is $\CC$\textit{-local},
if the natural map $\CC\to A/\m$ is an isomorphism. 
 
\section{Convergent NC power series}

\subsection{Absolutely convergent NC power series}

Refer to the Appendix for the notions of indices,
maps between indices,
noncommutative power series algebras, their morphisms  and
abelianization.
For each point $p\in \CC^n$, consider $x_i-p_i$ as a formal variable
and set 
$$\CC[[x-p]]:=\CC[[x_1-p_1\bbb x_n-p_n]].$$

A noncommutative power series 
$f=\sum_Ia_I (x-p)^I$ in $\CC[[x-p]]$
is called \textbf{absolutely convergent}, 
if the abelianization $f_{\abs}$ of $\sum_I|a_I| (x-p)^I$
is convergent. Write $$\CC\{x-p\}^{\abs}$$
for the local subalgebra of $\CC[[x-p]]$ 
of absolutely convergent noncommutative power series, and
for $r>0$, write
$$\CC\{x-p\}_r^{\abs}$$ for the subalgebra of $\CC\{x-p\}^{\abs}$ of those
$f$ such that $f_{\abs}$ converges on the open polydisk
$P(p,r)$ in $\CC^n$ with poly-radius $(r\ddd r)$.
A commutative power series $g=\sum_{\nu}b_\nu x^{\nu}$
converges on the open polydisk $P(r)=P(0,r)$, if and only if 
the series $g_{\abs}:=\sum_{\nu}|a_{\nu}|x^\nu$ converges on
$P(r)$. If a noncommutative power series $f=\sum_Ia_I x^I$ is
absolutely convergent on $P(r)$ and if
$|b_I|\leq |a_I|$, for each multi-index $I$, then
the series $\sum_I b_I(x-p)^I$ converges absolutely on $P(r)$.

\begin{observe}
If $g$ and each $f_1\ddd f_n$ are absolutely convergent power series, 
$g(f^1\ddd f^n)$ is absolutely convergent. 
\end{observe}
\begin{proof}
This follows by 
$g(f^1\ddd f^n)_{\ab}=g_{\ab}(f^1_{\ab}\ddd f^n_{\ab})$ and
equation~(\ref{einsetz}).
\end{proof}

\subsection{The $\MM$-spectrum}

Consider $\CC^k$ as normed space with respect to the maximum
norm and equip the 
matrix algebras $\MM=M_k(\CC)$ over $\CC$ with
the corresponding operator norm $||\cdot ||$.\\

Let $\C$ be a category of topological 
algebras. For objects $A$ and $\MM$, we call
$$\spec^{\MM}(A):=\Hom_{\C}(A,\MM)$$
the $\MM$\textbf{-spectrum} of $A$. For $\MM=M_k(\CC)$, the $\MM$-spectrum is
just the set 
$\operatorname{rep}_k(A)$ of continuous $k$-dimensional 
representations of $A$.
To each element $a\in A$, we
assign the canonical map 
$$\tilde a:\spec^{\MM}(A)\to \MM.$$
The ideal
$$I_A(\MM):=\{a\in A:\;\tilde a\equiv 0\}$$
is called \textbf{ideal of identities for} $\MM$.
For $\MM=M_k(\CC)$, we equip $\spec^{\MM}(A)$ with the weakest
topolology making the maps 
\begin{align*}
\tilde a:\spec^{\MM}(A)&\to\MM\\
\chi&\mapsto \chi(a)
\end{align*}
continuous, for each $a\in A$. 

\subsection{The joint spectral radius}

For an $n$-tuple $M=(M_1\ddd M_n)$ in
$\MM^n$ the value
$$|M|:=\limsup_{k\to\infty}\{\sqrt[|I|]{||M^I||}:\;|I|\leq k\}$$
is called \textbf{joint spectral radius} of $M_1\ddd M_n$
(see \cite{Theys}).
For elements $Q\in\CC^n$, we have
$|Q|=\max_{i=1\ddd n}|Q_i|$.
Observe that, for $k>1$, the
function $|\cdot |$ is not a semi-norm on $M_k(\CC)^n$. 
Let, for example,
$E_{ij}$ be the elementary matrices in $M_2(\CC)$, then
$|E_{12}+E_{12}|=1>0=|E_{12}|+|E_{21}|.$ 

\begin{satz}\label{theys}
The joint spectral radius $|\cdot|:\MM^n\to\RR_{\geq 0}$
has the following properties:
\begin{enumerate}
\item
It is continuous.
\item
$|M|=\lim_{k\to\infty}\max_{|I|=k}\sqrt[k]{||M^I||}$.
\end{enumerate}
\end{satz}
\begin{proof}
Both statements are proved in \cite{Theys}. The continuity was first
shown by \cite{Bara}.
\end{proof}

\begin{prop}\label{dreieck}
For elements $Q$ in $\CC^n$ and
$M$ in $\MM^n$, we
have the inequality
$$|M+Q|\leq m+|Q|,$$
where $m:=\max_{i=1\ddd n}||M_i||$.
\end{prop}
\begin{proof}
Set $q:=|Q|$.
Since for central elements $y_1\ddd y_n$,
\begin{equation}\label{binom}
(x+y)^I=\sum_{\alpha:J\to I} x^J y^{I-_\alpha J},
\end{equation}
we have
\begin{equation}\label{feins}
(x+\eins)^I=\sum_{J\leq I}{ I\choose J} x^J,
\end{equation}
where $\eins=(1\ddd 1)$.
By equation (\ref{binom}),
\begin{align*}
||(M+Q)^I||=&||\sum_{\alpha:J\to I} M^J Q^{I-_{\alpha}J}||\\
\leq&\sum_{\alpha:J\to I}||M^J||\cdot q^{|I|-|J|}\\
=&\sum_{J\leq I}{I \choose J} ||M^J||\cdot q^{|I|-|J|}\\
\leq&\sum_{J\leq I}{I \choose J} ||M||^J\cdot q^{|I|-|J|},
\end{align*}
where $||M||:=(||M_1||\ddd ||M_n||)$. 
Thus,
\begin{equation*}
\sqrt[|I|]{||(M+Q)^{I}||}\leq
q\cdot\sqrt[|I|]{\sum_{J\leq I}{I\choose J} ||\frac{M}{q}||^J}
\end{equation*}
By equation (\ref{feins}), the last expression equals
\begin{equation*}
q\cdot\sqrt[|I|]{(||\frac{M}{q}||+\eins)^I}\leq 
q(\frac{m}{q}+1)
=m+q.
\end{equation*}
\end{proof}

For $P\in\MM^n$ and $r>0$, the subset
$$P^{\MM}(P,r):=\{M\in\MM^n:\; |M-P|<r\}$$
of $\MM^n$ will be called \textbf{stable region in} $\MM^n$ 
with radius $r$ and center $P$. By Theorem~\ref{theys},
stable regions are open in $\MM^n$. Like open polydisks in $\CC^n$,
they play the role of smooth affine spaces.
There are algorithms \cite{BlNe} for arbitrary good approximations
of the joint spectral radius of n-tuples $M$ of matrices
with rational entries.
It is not known, so far, if the question 
weather or not an element $M\in\MM^n$ 
satisfies $|M|<1$, is decidable.
However, there is the following negative result
\cite{TsBl}:

\begin{satz}[Tsitsikilis-Blondel]
It is undecidable, if $|M|\leq 1$,
for a given n-tuple $M$ of square
matrices with rational entries.
\end{satz}

\subsection{Convergence}

For $\MM=M_k(\CC)$, we fix the basis
$E_1\ddd E_{k^2}$, where $E_1$ is the identity matrix
and $E_2\ddd E_{k^2}$ are the lexicographically ordered
elementary matrices $E_{ij}$ with $(i,j)\neq(1,1)$.\\

The identification  $\tau:\MM^n\to\CC^{k^2\cdot n}$ is
a map of topological vector spaces. 
Consider the algebra $\CC\{x_{11}\ddd x_{nk^2}\}_{\MM,r}$
of all commutative power series in $x_{11}\ddd x_{nk^2}$,
absolutely convergent at each point of $\tau(P^{\MM}(r))$.
The algebra $\CC\{x_{11}\ddd x_{nk^2}\}$ is Fr\^echet with
respect to the topology of compact convergence.
The elements of the algebra
$$\MM[[x_{11}\ddd x_{nk^2}]]:=
\CC[[x_{11}\ddd x_{nk^2}]]\ot_{\CC}\MM$$ 
correspond to $k^2$-tuples of elements in
$\CC[[x_{11}\ddd x_{nk^2}]]$
via the coordinate projections $\MM\to\CC$. The
algebra
$$\MM\{x_{11}\ddd x_{nk^2}\}_r:=
\CC\{x_{11}\ddd x_{nk^2}\}_{\MM,r}\ot_{\CC}\MM$$ 
is again Fr\^echet. 
It plays the role of a universal object,
in the sense that any $k$-dimensional
representation of the ``convergent power series algebra''
(see below) factors through it.\\

Consider the map
$$\phi^{\MM}:\CC[[x_1|\ldots|x_n]]\to
\MM[[x_{11}\ddd x_{nk^2}]],$$
sending $x_i$ to the ``generic matrix'' 
$\sum_{j=1}^{k^2}E_j\cdot x_{ij}$.
Set
$$\CC\{x\}^{\MM}_r:=
(\phi^{\MM})^{-1}(\MM\{x_{11}\ddd x_{nk^2}\}_r).$$ 
and
$$\CC\{x\}_r:=\cap_{\MM}\CC\{x\}^{\MM}_r.$$
The elements of
$$\CC\{x\}:=\cup_{r>0}\CC\{x\}_r.$$
will be called \textbf{convergent} noncommutative power
series.\\ 

Observe that, for $M\in P^{\MM}(r):=P^{\MM}(r,0)$ and $f=\sum_I a_Ix^I$
in $\CC\{x\}^{\abs}_r$,
the sum $\sum_Ia_IM^I$ converges absolutely in $\MM$.
Thus, each absolutely convergent noncommutative
power series is convergent. 
By Example~\ref{stand} below,
the converse is not true:

\subsection{Ideals of matrix identities}

The kernel of the 
restriction of $\phi^{\MM}$
to the noncommutative polynomial ring
$\CC[x]$ is the ideal of matrix identities
$I_{\CC[x]}(\MM)$, that was studied
by Procesi \cite{Procesi} and others. 
Since $\phi^{\MM}$ respects total degrees, an NC polynomial
$f$ is in $I_{\CC[x]}(\MM)$, if and only if each of its homogeneous
components $f_d$ of total degree d is in $I_{\CC[x]}(\MM)$.
By the Amitsur-Levitzki Theorem \cite{AmiLe} (see 
Example~\ref{stand} below),
 there is an identity of degree $2k$,
for $\MM=M_k(\CC)$.
It is known \cite{Procesi} that 
$I_{\CC[x]}(\MM)$
contains no homogeneous polynomial of degree $\leq 2k$.\\

For $A\in\{\CC\{x\}_r,A=\CC\{x\}^{\abs}_r,\CC[x]\}$, we just write
$I(\MM)$ for 
the ideal $I_A(\MM)$ of matrix identities, or just $I_k$ if $\MM=M_k(\CC)$.
Again,
$$I(\MM)=\Kern(\phi^{\MM})|_A.$$

\begin{prop}\label{vergl}
$I_k\sub (x)^{2k}$
\end{prop}
\begin{proof}
For $A=\CC[x]$, the statement is know to be true.
By homogenity of the maps $\phi^{\MM}$, 
a power series $f$ with homogeneous components $f_d$
of total degree $d$ belongs to $I_k$, if and only if
each $f_d$ belongs to $d$. Thus if $f\in I_k$, each $f_d$ with
$d<2k$ is trivial.
\end{proof}

\begin{kor}
$\cap_k I_k=0.$
\end{kor}

\begin{beisp}\label{stand}
For $l\geq 1$, let $s_l\in\ZZ[x_1|\ldots|x_l]$ be the standard identity 
$$s_l:=\sum_{\sig\in S_l}\sgn(\sigma)x_{\sig(1)}\ppp x_{\sig(l)}.$$
One can verify as an exercise that
each $l-1$-dimensional $\CC$-algebra satisfies
\begin{equation}\label{standard}
\forall x_1\ddd x_l:\;s_l(x_1\ddd x_l)=0.
\end{equation}
The Amitsur-Levitzki Theorem \cite{AmiLe}
states that $M_k(\CC)$
satisfies (\ref{standard}), for $l=2k$.
The algebra map
\begin{align*}
\CC[z_1|\ldots|z_l]&\to\CC[x|y]\\
z_i&\to x y^i
\end{align*}
is injective (see Lam).
In consequence, $f_l:=s_l(xy\ddd xy^l)$ is a non-trivial identity 
of degree $\frac{l(l+1)}{2}$
for all $\CC$-algebras of dimension $\leq l-1$ (and for $M_{l/2}(\CC)$,
if $l$ is even). We have that
$$f_l=\sum_{|I|}a_I(x,y)^I,$$
where $a_I\in\{0,\pm 1\}$ is non-trivial for exactly
$l!$ many indices. The series
$$f:=\sum_{l\geq 1}f_l$$
is not absolutely convergent, since $f_{\abs}=\sum_ll!x^ly^{l(l-1)/2}$.
It is convergent, since $\phi^{\MM}(f)$ is a polynomial, for
each $\MM=M_k(\CC)$.
\end{beisp}

\subsection{The topology of compact convergence}

For
compact subsets $K\sub P^{\MM}(r)$, we define a
seminorm $||\cdot ||_K$ on $\CC\{x\}_r$ via
$$||f||:=||\phi^{\MM}(f)||_K.$$
The induced topology is already defined by a countable 
subfamily of these seminorms and is called \textbf{topology
of compact convergence}.
By the following
theorem, it is a Fr\^echet-topology.

\begin{lemma}\label{koeffi}
If a sequence $(f_m)_m$ in $\CC\{x\}_r$ with $f_m=\sum_Ia_{m,I}x^I$
is Cauchy with respect to the topology of compact convergence, then,
for each $I$, 
the coefficient sequence $(a_{m,I})_m$ converges in $\CC$.
\end{lemma}
\begin{proof}
The corresponding statements for commutative convergent
power series is known to be true.
Since $I(\MM)\sub(x)^{2k}$, for $\MM=M_k(\CC)$, the restriction 
of $\phi^{\MM}$ defines a monomorphism 
\begin{equation}\label{injek}
\CC[x_1|...|x_n]_{\leq 2k-1}\to\MM[x_{11}\ddd x_{nk^2}]_{\leq 2k-1}
\end{equation}
of finite dimensional vectorspaces.
By assumption,
$(\phi^{\MM}(f_m))_m$ is a Cauchy sequence in
$\MM\{x_{11}\ddd x_{nk^2}\}_r$, thus it converges coefficientwise,
since $\phi^{\MM}(f_m)$ is just a matrix of commutative convergent
powers series.
Injectivity of the morphism (\ref{injek}) implies that
the series $(a_{m,I})_m$ converges, for each $|I|\leq 2k-1$.
\end{proof}

\begin{satz}\label{frechet}
The algebra $\CC\{r\}_r$ is complete with respect to the topology
of compact convergence.
\end{satz}
\begin{proof}
Consider a Cauchy sequence $f_k$ in $\CC\{x\}_r$ with
$f_k=\sum_I a_{k,I}x^I$.
By Lemma~\ref{koeffi}, 
for each multi-index $I$, the sequence 
$(a_{k,I})$ has a limit $a_I\in\CC$. 
Set $f:=\sum_Ia_Ix^I$.
Since $f_k$ converges coefficientwise to $f$, 
the Cauchy sequence $\phi^{\MM}(f_k)$ converges coefficientwise
to $\phi^{\MM}(f)$. By completeness, the series 
$\phi^{\MM}(f)$ belongs to $\MM\{x_{11}\ddd x_{nk^2}\}_r$.
\end{proof}

\section{Representations}

\subsection{$\MM$-spectra of the convergent power series algebra}

Consider the induced topology on the subalgebra
$\CC\{x\}_r^{\abs}\sub\CC\{x\}_r$ of absolutely convergent power series.
Recall that the $\MM$-spectrum $\spec^{\MM}(\CC\{x\}^{\abs}_r)$
is the set of all continuous algebra homomorphisms
$f:\CC\{x\}_r^{\abs}\to\MM$. Any such $f$ is uniquely defined
by the values $f(x_1)\ddd f(x_n)$ on the variables.

\begin{satz}\label{homeom}
The natural map $P^{\MM}(r)\to\spec^{\MM}(\CC\{x\}^{\abs}_r)$
ist a homeomorphism.
\end{satz}
\begin{proof}
To prove bijectivity,
it suffices to show that if $\M\in\MM^n$ is not in $P^{\MM}(r)$, there
is an element $f\in\CC\{x\}^{\abs}_r$ such that the sequence
$$(||\sum_{|I|=k}a_I M^I||)_{k\geq 0}$$ 
does not converge to zero.
If $|M|\geq r$, we may assume that there exists 
a sequence $(I(m))_m$ in $G$
with $|I(m+1)|>|I(m)|$, for each $m\geq 0$, such that the sequence
$(\sqrt[|I(m)|]{||M^{I(m)}||})_m$ is monotouneously increasing and
converges to $r$. The one variable
power series 
$$\sum_m\frac{1}{||M^{I(m)}||}z^{|I(m)|}$$
converges for $0<z<r$: Assume that $z<r-\delta$, for a
$\delta\in(0,r)$. For $m>>0$, we have that 
$||M^{I(m)}||^{-1}<(\frac{1}{r-\delta})^{|I(m)|}$,
thus 
$$\frac{z^{|I(m)|}}{||M^{I(m)}||} <(\frac{z}{r-\delta})^{|I(m)|}.$$
In consequence, if we set $a_I:=\frac{1}{||M^{I(m)}||}$, for
$I=I(m)$ and $a_I=0$, if $I$ is not in $\{I(m):\;m\geq 0\}$, the
series
$\sum_Ia_Ix^I$ belongs to $\CC\{x\}^{\abs}_r$. On the other hand-side,
the sequence $$(||\sum_{|I|=k}a_I M^I||)_{k\geq 0}$$ 
does not converge to zero.
It is left to the reader to verify that the map is bicontinuous.
\end{proof}

Observe that
the inclusion $\CC\{x\}^{\abs}_r\to\CC\{x\}_r$ induces an isomorphism
on each $\MM$-spectrum.

\begin{kor}
For $r>0$ and $\MM=M_k(\CC)$,
the topological space $\spec^{\MM}(\CC\{x\}_r)$ of continuous
$k$-dimensional representations of $\CC\{x\}_r$ is just
the stable region $P^{\MM}(r)$ in $\MM^n$.
\end{kor}

\subsection{The $I_\bullet$-topology}

For $A\in\{\CC\{x\}_r,A=\CC\{x\}^{\abs}_r,\CC[x]\}$,
the filtration
$$ \ldots I_3\sub I_2\sub I_1=[A,A]\sub A$$
defines the $I_\bullet$-\textbf{topology} on $A$.

\begin{satz}\label{volls}
The algebra $\CC\{x\}_r$ is complete with respect to
the $I_\bullet$-topology.
\end{satz}
\begin{proof}
Consider an $I_\bullet$-Cauchy sequence $(f_k)$ in $\CC\{x\}_r$.
Denote $f_k^{(m)}$ the homogeneous part of $f_k$ of total degree $m$.
Replacing the sequence by an appropriate subsequence, we may
assume that $f_k-f_j\in I_{k+1}\sub (x)^{2(k+1)}$ for $j\geq k$.
This implies the following identities:
\begin{enumerate}
\item\label{kr1}
$\phi^{\MM}(f_k)=\phi^{\MM}(f_j)$, for $\MM=M_k(\CC)$ and $j\geq k$.
\item\label{kr2}
$f_j^{(j)}=f_k^{(j)}$, for $j\leq k$.
\end{enumerate}
  Set $f^{(k)}:=f^{(k)}_{k}$, for $k\geq 0$ and
$f:=\sum_{k\geq 0}f^{(k)}$. 
First we show that $f-f_k\in I_k$, for each $k\geq 0$, and in
consequence, that the series $(f_k)$ converges to
$f$ in the $I_\bullet$-topology of $\CC[[x]]$.
By equation (\ref{kr1}), for $\MM=M_k(\CC)$ and $j>k$,
we have that
\begin{equation*}
\phi^{\MM}(f-f_k)^{(j)}=
\phi^{\MM}(f^{(j)}_j-f^{(j)}_k)=
\phi^{\MM}(f_j-f_k)^{(j)}
=0.
\end{equation*}
By equation (\ref{kr2}), 
for $\MM=M_k(\CC)$ and $j\leq k$.
we have that
\begin{equation*}
\phi^{\MM}(f-f_k)^{(j)}=
\phi^{\MM}(f^{(j)}_j-f^{(j)}_k)=0.
\end{equation*}
In consequence, $\phi^{\MM}(f-f_k)=0$.
To conclude the proof, we must show that $f$ belongs
to $\CC\{x\}_r^{\MM}$, for each $\MM=M_l(\CC)$:
For arbitrary $P\in P^{\MM}(r)$ and $\epsilon>0$, we can find a
$k_0\geq l$ such that, for each $k_1>k_o$, we have that
$$\sum_{k=k_0}^{k_1}\phi^{\MM}(f_l)^{(k)}_{\abs}(P)\leq\epsilon.$$
Since
\begin{align*}
\phi^{\MM}(f_l)^{(k)}&=\phi(\MM)(f_k)^{(k)}\quad\text{ (by equation 
(\ref{kr1}))}\\
&=\phi^{\MM}(f_k^{(k)})\\
&=\phi^{\MM}(f^{(k)})=\phi^{\MM}(f)^{(k)},
\end{align*}
we have that
$$\sum_{k=k_0}^{k_1}\phi^{\MM}(f)^{(k)}_{\abs}(P)\leq\epsilon.$$
This finishes the proof.
\end{proof}

Let $\sC$ be the category all Fr\^echet algebras of the form
$\CC\{x\}_r/J$, where $J$ is a Fr\^echet-closed ideal.

\begin{kor}
$$\CC\{x\}_r={\lim_k}^{\sC}\CC\{x\}_r/I_k.$$
\end{kor}

\section{Appendix: Formal NC power series}

\subsection{Multi-indices}

Fix a dimension $n\geq 0$.
Let $G$ be the semi-group of 
\textbf{(noncommutive) multi-indices,} i.e. of tuples 
$I=(i_1\ddd i_{|I|})$ with $i_k\in\{1\ddd n\}$
of length $|I|\geq 0$. The addition of 
multi-indices $I,J$ is defined by $I+J:=(i_1\ddd i_{\#I},j_1\ddd j_{\#J})$.
The zero element is the empty multi-index denoted by $0$. 

Consider the partial relation on $G$ defined by
$I<J$ iff and only if
$|I|<|J|$ and $I$ is obtained from $J$ by deletion of indices.
If $I\leq J$ (i.e. $I<J$ or $I=J$), there exists
an order-preserving, injective map $\alpha$ from
the set $\{1\ddd\#I\}$ into $\{1\ddd\#J\}$ such that
$i_k=j_{\alpha(k)}$, for $k=1\ddd\#I$. 
We say that $\alpha$ is a map from $I$ to $J$.
Set ${J \choose I}$ to be the number of such maps $\alpha:I\nach J$. 
By definition, ${J \choose 0}=1$, for any $J$.
We have ${J \choose I}\leq {|J|\choose|I|}$.
Denote $J-_\alpha I$ be the multi-index obtained from $J$ by deletion
of $j_{\alpha(1)}\ddd j_{\alpha(\#I)}$.
Let $G_{\ab}$ denote the commutative semi-group of $n$-tuples 
$\nu=(\nu_1\ddd\nu_n)$ in $\{0,1,2,\ldots\}$.
We have an epimorphism $\ab:G\to G_{\ab}$ of semigroups
sending $I\in G$ to the multi-index $\nu\in G_{\ab}$ with
$\nu_k=\{l:i_l=k\}$.

\subsection{The category of formal NC power series}

For $n>0$, 
write $\CC[[x_1\bbb x_n]]$ or just $\CC[[x]]$
for the noncommutative
power series algebra in formal variables $x_i$.
For each multi-index $I\in G$, set
\begin{align*}
x^I&:=x_{i_1}\ppp x_{i_k}\quad\text{ and }\\
x^0&:=1.
\end{align*}

If $g=\sum_Jb_J x^J$ is a power series in $\CC[[x]]$
and if the power series
$f^1\ddd f^n$ with $f^k=\sum_Ia_I^k y^I$
belong to the maximal ideal $(y)$
of
$\CC[[y]]=\CC[[y_1\bbb y_m]]$, we can form the power series
\begin{equation}\label{einsetz}
g(f^1,...,f^n):=\sum_J(\sum b_K\cdot
a^{k_1}_{I_{k_1}}\ppp a^{k_{|K|}}_{I_{k_{|K|}}})y^J,
\end{equation}
where the sum in the bracket is over all multi-indices
$K,I_1\ddd I_{|K|}$ such that $J=I_{k_1}\kkk I_{k_{|K|}}$.
A \textbf{morphism} $\CC[[x]]\nach\CC[[y]]$
of noncommutative power series algebras is a local algebra
homomorphism of the form $g\mapsto g(f_1\ddd f_n)$,
for a given $n$-tuple $(f^1\ddd f^n)$ of elements of
$(y)$.

\begin{observe}\label{composition}
Let $f:\CC[[x]]\nach\CC[[y]]$ and $g:\CC[[y]]\nach\CC[[z]]$
be given by $f(x_s)=\sum_Ia^s_I y^I$ and
$g(y_t)=\sum_{J}b^t_J z^J$. Then, 
$g\circ f(x_s)=\sum_Kc_K z^K$, with
$$c_K=\sum_Ia^s_I\sum b^{i_1}_{J_1}\ppp b^{i_{|I|}}_{J_{|I|}},$$
where the second sum is taken over all multi-indices
$J_1\ddd J_{|I|}$ such that $J_1\kkk J_{|I|}=K$.
\end{observe}

For an endomorphism $f$ of
$\CC[[x]]$, with $f(x_k)=\sum_Ia^k_{I}x^I$, set
$Jf$ to be the $n\times n$-matrix $(a^k_{(i)})_{k,i}$.
\begin{satz}
An endomorphism $f$ of $\CC[[x]]$ is an automorphism, if and
only if $Jf$ is invertible.
\end{satz}
\begin{proof}
Suppose that $Jf$ is invertible. Inductively, for
$k\geq 1$, $l=1\ddd k$ and $J\in\{1\ddd n\}^l$,
we  will define coefficients $b^l_{J}\in\CC$
such that the endomorphism $g^{k}$ of $\CC[[x]]$ with
$g^{k}(x_l)=\sum_{|J|\leq k}b^l_{J}x^J$ is inverse to $f$
modulo $(x)^{k+1}$.
For $k=1$, by Observation~\ref{composition}, the necessary
(and sufficient) conditions are
\begin{align*}
\sum_{i=1}^na^s_{(i)}b^i_{(s)}=&1\quad\text{ for all } s,\\
\sum_{i=1}^na^s_{(i)}b^i_{(s)}=&0\quad\text{ for all } s \neq l.
\end{align*}
We can find such coefficients, if and only if $Jf$ is invertible.
Now suppose that $Jf$ is invertible and that the coefficients
$b^l_{J}$ are constructed adequately, for $|J|\leq k-1$.
For fixed $J$ with $|J|=k$, we have to
find $b^l_{J}$ such that 
$$\sum_{i=1}^n a^s_{(i)}b^i_{J}+\sum_{|I|\geq 2}a^l_{I}
\sum b^{i_1}_{J_1}\ppp b^{i_{|I|}}_{J_{|I|}}=0.$$
The second sum is known, and since $Jf$ is invertible, we
can find adequate $b^1_{J}\ddd b^n_{J}$.
\end{proof}

We have a canonical epimorphism $\ab$
from $\CC[[x]]=\CC[[x_1\bbb x_n]]$
to the commutative power series algebra 
$\CC[[x]]_{\ab}=\CC[[x_1\ddd x_n]]$.
We write  $f_{\ab}$ instead of $\ab(f)$.

\begin{kor}
Each lift of an automorphism of a commutative formal
power series algebra to an endomorphism of the non-commutative
formal power series algebra  is again an automorphism. 
\end{kor}

\subsection{Products}
For power series algebras $\CC[[x]]=\CC[[x_1|...|x_n]]$ and 
$\CC[[y]]=\CC[[y_1|...|y_m]]$, we
define the \textbf{free product} 
$$\CC[[x]]\ast \CC[[y]]:=\CC[[x_1|...|x_n|y_1|...|y_n]]$$
and the \textbf{complete tensor product} 
$\CC[[x]]\otc_\CC\CC[[y]]$ as the power series
algebra in $x_1,...,x_n$ and $y_1,...,y_m$, where the $y_j$ are 
assumed to commute with the $x_i$.

\subsection{Finitely generated ideals}

For noncommutative power series algebras, the concept of finitely generated
two-sided ideals has to be slightly adapted. As a reason,
we give the following example:
\begin{beisp}
Let $(x)$ be the two-sided ideal of $\CC[[x|y]]$,
consisting of all noncommutative power series
where at least one factor $x$ arises in each monomial.
Observe that $(x)$ is not the two-sided ideal generated by $x$ in
the algebraic sense.
\end{beisp} 
\begin{proof}
Assume that $(x)$ is the two-sided ideal generated by $x$ in the 
algebraic sense. Then we can find power series 
$f_i=\sum_j a_{ij}y^j$ and $g_i=\sum_jb_{ij}y^j$ in
$\CC[[y]]$, $i=1\ddd N$ such that
$$xyx+y^2xy^2+...=\sum_{i=1}^Nf_i(y)\cdot x\cdot g_i(y).$$
The right hand-side takes the form
$\sum_{j,k}(\sum_i a_{ij}b_{ik})y^jxy^k$.
In particular, for $j,k\leq N+1$, we would get
$\sum_i^Na_{ij}b_{ik}=\delta_{j,k}$,
which is impossible, since the left hand-side is
a product of two matrices of rank at most $N$.
\end{proof}

By definition, the opposite algebra 
$A^{\op}$ of an algebra $A$
is the set $\{a^{\op}:\;a\in\CC[[x]]\}$ with operations
$a^{\op}_1+a^{\op}_2=(a_1+a_2)^{\op}$ and 
$a_1^{\op}\cdot a_2^{\op}=(a_2\cdot a_1)^{\op}$. Observe that
$\CC[[x]]^{\op}$ is naturally isomorphic to the power
series algebra $\CC[[x_1^{\op}|...|x_n^{\op}]]$ and the assignment
\newcommand{\OP}{\operatorname{OP}}
$x\mapsto x^{\op}$ defines an isomorphism 
$\OP:\CC[[x]]\nach \CC[[x]]^{\op}$.
Attention, in general, $\OP(f)\neq f^{\op}$, for example, 
$\OP(x_1 x_2)=x_1^{\op} x_2^{\op}=(x_2 x_1)^{\op}$.
We define the \textbf{(complete) envelopping algebra} 
$\CC[[x]]^{\ec}$ of $\CC[[x]]$
as $\CC[[x]]\otc_\CC\CC[[x]]^{\op}$.
Consider the natural epimorphism
$$\alphac:\CC[[x]]^{\ec}\nach \CC[[x]].$$
For a two-sided ideal 
$J\sub \CC[[x]]$, the inverse image $\alphac^{-1}(J)$ is not, in general
a left ideal of $\CC[[x]]^{\ec}$, since it is not, in general, 
closed under left multiplication by elements of $\CC[[x]]^{\ec}$.
We define the \textbf{completion} $\hat{J}$ of $J$ as the image under
$\alphac$ of the $\CC[[x]]^{\ec}$-left ideal generated by 
$\alphac^{-1}(J)$.
For simplicity, for elements $f_1,...,f_m$ in $\CC[[x]]$, we shall write
$(f_1,...,f_m)$ for the completion of the
two-sided ideal generated by the $f_i$. A two-sided ideal of
the form $(f_1,...,f_m)$ will be called \textbf{finitely generated}.

\begin{prop}
The Kernel $\K$ of the abelization 
$\ab:\CC[[x_1|...|x_n]]\nach\CC[[x_1\ddd x_n]]$ is
finitely generated by the commutators $[x_i,x_j]=x_ix_j-x_jx_i$,
for $1\leq i<j\leq n$.
\end{prop}
\begin{proof}
We show that for each noncommutative power series $f$,
the difference $f-f_{\ab,\unab}$ is in the image
under $\alphac$ of the $\CC[x]^{\ec}$-left ideal
generated by the commutators $[x_i,x_j]$, for $i<j$.
Without restriction, let
$f$ be the sum of its homogeneous components
$f_k$ of degree $k\geq 2$.
Each difference $f_k-f_{\ab,\unab,k}$ is of the
form $\alphac(\sum_{i<j}c^{ij}_k\cdot [x_i,x_j])$, for
certain homogeneous $c_k^{i,j}$ in $\CC[x]^{\ec}$ of
degree $k-2$.
Thus $f-f_{\ab,\unab}$ is the image under $\alphac$ of
$\sum_{i<j}(\sum_k c_k^{ij})\cdot[x_i,x_j]$.  
\end{proof}

\subsection{Locality}

A family $(h_\alpha)_{\alpha\in A}$ of power series in $\CC[[x]]$
is called \textbf{summable}, if, for each multi-index
$I$, there are only finitely many $\alpha\in A$ such that
$h_{\alpha,I}\neq 0$. In this case, we can form the power
series
$$\sum_{\alpha\in A}h_\alpha:=\sum_I(\sum_\alpha h_{\alpha,I})(x-p)^I.$$

\begin{prop}\label{localring}
The algebra $\CC[[x]]$  is local with maximal
ideal $(x)$ generated by $x_1\ddd x_n$.
\end{prop}
\begin{proof}
It suffices to show that each element $f$ of $\CC[[x]]\setminus (x)$
is a unit. Without restriction, say $f_0=1$.
Then the family $(1-f)^j;j\geq 0$ is summable.
We have 
\begin{align*}
f\cdot\sum_{j=0}^\infty(1-f)^j&=(1-(1-f))\sum_{j=0}^\infty(1-f)^j=\\
&=\sum_{j=0}^\infty(1-f)^j-\sum_{j=1}^\infty(1-f)^j=1.
\end{align*}
Thus $f$ is a left unit. 
In the same way, we show that $f$ is a right unit.
If $f$ is convergent, the sum $\sum(1-f)^j$ is also convergent. This
follows exactly as in the commutative case.
\end{proof}


\begin{thebibliography}{lalalalala}


\bibitem{AmiLe} S.A. Amitsur; J. Levitzki: \textit{Minimal identities for
algebras,} Proc. AMS 1, 449-463 (1950).

\bibitem{Bara} N.E. Barabanov: \textit{Lyapunov inducator of discrete
inclusions I-III,} Autom. Remote Control 49, 2: 152-157,
3: 283-287, 5: 558-565 (1988).

\bibitem{BlNe} Vincent D. Blondel; Yurii Nesterov: \textit{Fast and 
precise approximations of the joint spectral radius,} CORE discussion paper
(2003), available at www.core.ucl.ac.be/services/COREdp03.html.
 

\bibitem{Cohn} P.M. Cohn: \textit{Free rings and their relations,}
Academic Press (1971).

\bibitem{Grau} Hans Grauert; Reinhold Remmert:
\textit{Theorie der Steinschen R\"aume,}
Springer  (1977).


\bibitem{Kaly} Dmitry S. Kalyuzhny\u{i}-Verbovetzki\u{i}:
\textit{Carath\'eodory interpolation on the non-commutative
polydisk,} arxiv:math.FA/0412161v2 (2005).

\bibitem{Kapr} Mikhail Kapranov:
\textit{Noncommutative geometry based on commutator expansions,}
J. Reine Angew. Math. 505, 73-118 (1998).

\bibitem{Lam} T.Y. Lam: \textit{A first course in noncommutative rings,}
2nd edition, Springer (2001).

\bibitem{LeBruyn} Lieven Le Bruyn: \textit{Noncommutative geometry@n,}
available at www.math.ua.ac.be/~lebruyn.

\bibitem{Procesi} Claudio Procesi: \textit{Rings with polynomial identities,}
Dekker (1973).

\bibitem{ncaspaces} Frank Schuhmacher: \textit{Noncommutative complex
analytic spaces,} math/QA.0606150.


\bibitem{spaces} Frank Schuhmacher: \textit{Noncommutative analytic spaces,}
to appear.

\bibitem{Theys} Jacques Theys: \textit{Joint spectral radius: 
theory and approximations,} Phd Thesis, Universit\'e Catholique de Louvain
(2005).

\bibitem{TsBl}
 J.N. Tsitsiklis; V.D. Blondel: The Lyapunov exponent and joint spectral
radius of matrices are hard - when not impossible - to compute and to
approximate, Math. of Control, Signals and Systems 10, 31-40 (1997).

\end{thebibliography}
\end{document}